\documentclass{ifacconf}

\usepackage{graphicx}      
\usepackage{natbib}        
\usepackage{amsmath}
\usepackage{amssymb}
\usepackage{xcolor}
\usepackage{url}
\usepackage{natbib}

\newcommand{\mb}[1]{\mathbf{#1}}
\newcommand{\vo}[1]{\mathbf{#1}}
\newcommand{\x}{\mb{x}}

\newcommand{\xdot}{\dot{\mb{x}}}

\newcommand{\X}{\mb{X}}

\newcommand{\m}{\mb{m}}

\newcommand{\Sd}{\mb{S}}
\newcommand{\tSd}{\mb{\tilde{S}}}

\newcommand{\f}{\mb{f}}

\newcommand{\E}[1]{\mathbb{E}[#1]}

\newcommand{\nx}{d}
\newcommand{\nb}{n_{\set{B}}}
\newcommand{\nd}{n_{\set{D}}}
\newcommand{\Real}{\mathbb{R}}
\newcommand{\real}{\mathbb{R}}
\newcommand{\maxent}{\textit{maxent }}
\newcommand{\cvxh}[1]{\text{Conv}(#1)}
\newcommand{\set}[1]{\mathcal{#1}}

\newcommand{\param}{\vo{\Delta}}
\newcommand{\tparam}{\tilde{\param}}
\newcommand{\domainD}{\set{D}_{\param}}
\newcommand{\pdfp}{p(\param)}
\newcommand{\xpc}{\x_{\text{pc}}}

\newcommand{\xpcdot}{\xdot_{\text{pc}}}
\newcommand{\fhat}{\vo{\hat{f}}}

\newcommand{\phip}{\vo{\Phi}(\param)}
\newcommand{\phipt}{\vo{\Phi}^T(\param)}
\newcommand{\vect}[1]{\textbf{vec}\left(#1\right)}

\newcommand{\A}{\vo{A}}
\newcommand{\B}{\vo{B}}
\newcommand{\C}{\vo{C}}
\newcommand{\D}{\vo{D}}
\newcommand{\I}[1]{\vo{I}_{#1}}

\newcommand{\Ap}{\A(\param)}

\newcommand{\eqnlabel}[1]{\label{eqn:#1}}
\newcommand{\figlabel}[1]{\label{fig:#1}}
\newcommand{\eqn}[1]{(\ref{eqn:#1})}
\newcommand{\Eqn}[1]{Equation (\ref{eqn:#1})}
\newcommand{\fig}[1]{Fig.\ref{fig:#1}}

\begin{document}
\begin{frontmatter}

\title{Data-driven Solution of Stochastic Differential Equations Using Maximum Entropy Basis Functions\thanksref{footnoteinfo}}

\thanks[footnoteinfo]{This work was supported by the National Science Foundation (grant no. 1762825).}
\author{Vedang M. Deshpande \& Raktim Bhattacharya}
\address{Department of Aerospace Engineering, Texas A\&M University \\
   College Station, TX 77843, USA \\ (e-mail: vedang.deshpande@tamu.edu, raktim@tamu.edu)}

\begin{abstract}                
In this paper we present a data-driven approach for uncertainty propagation.
In particular, we consider stochastic differential equations with parametric uncertainty.
Solution of the differential equation is approximated using maximum entropy (maxent) basis functions similar to polynomial chaos expansions.
Maxent basis functions are derived from available data by maximization of information-theoretic entropy, therefore, there is no need to specify basis functions beforehand.
We compare the proposed maxent based approach with existing methods.
\end{abstract}


\end{frontmatter}

\section{Introduction}
Most physical systems are modeled by non-linear differential equations. Such models often suffer from parametric uncertainties, which in turn leads to inaccurate prediction of the states of the system.
Even if the state of a system is exactly known at a time instant, due to the parametric uncertainties, the states at subsequent time instants determined from uncertain models exhibit stochastic nature, and often one is interested in temporal evolution of distribution or statistical moments of the states.

In general, determining distribution of states, by propagation through stochastic models, is an infinite dimensional and computationally intractable problem. Polynomial chaos expansions (PCE) have been widely studied as an approach to construct a reduced order and computationally tractable surrogate model for the original system.
In PCE, a stochastic process is expressed as a weighted sum of polynomials of random variables. The polynomials or basis functions are selected according to the underlying distribution of the random parameters, and the optimal weights or coefficients of the chaos are determined using the so called Galerkin projection.
A correspondence between stochastic distributions and optimal basis from Askey family of orthogonal polynomials which leads to the exponential convergence of error with respect to the order of approximation is discussed in \cite{wiener_askey}. This is also known as generalized polynomial chaos (gPC).
Application of PCE for modeling uncertainties in physical systems has been discussed in many works, such as \cite{dXiu_jcp, sandia, fracture}, to name a few.

The classical polynomial chaos approach assumes that the distributions of involved random variables are exactly known. However, in reality, this assumption does not hold true. Often, limited information is available about the distributions in the form of finite number of samples, and engineering models must be developed using such limited set of data.
A data-driven approach for chaos expansion, termed as arbitrary polynomial chaos (aPC), is discussed in \cite{aPC}. In aPC, instead of fitting a probability distribution function over the given data, an orthogonal basis of polynomials for chaos expansion is constructed from raw moments of the available samples.
In the present work, we pursue a similar approach, but instead of constructing polynomials from raw moments, we use maximum entropy (\textit{maxent}) functions which are derived from the available samples, as basis to construct a chaos expansion.

The principle of maximum entropy was first introduced in \cite{jaynes1,jaynes2} as an approach to draw the least biased inference from incomplete or limited data.
Maximum entropy based methods are often used for inferring a probability distribution from sparse information and it has found applications in many diverse fields such as structural mechanics, image processing, and machine learning, among others. See \cite{maxent_struct, img_proc, maxent_irl}.
The \maxent basis functions that we use in this work arise as shape functions in the context of \maxent based polygonal interpolation, as discussed in \cite{sukumar,ortiz}. For the sake of completeness, \maxent based polygonal interpolation is briefly reviewed in Section \ref{maxent_fcn}.
Our recent work, \cite{maxentFunApprox}, discussed the data-driven function approximation using such \maxent basis functions in the context of dynamics modeling. The results show that \maxent based approach can model the unknown dynamics with good accuracy even if the available data set is sparse.
This served as a motivation to pursue the present study.

The objective of this work is to investigate the accuracy and convergence properties of chaos expansions constructed using \maxent basis functions for solution of stochastic differential equations (SDEs). 

The key highlights of this work are as follows.
\begin{enumerate}
\item[$\bullet$] We present a novel approach for developing surrogate models of SDEs using data-driven \maxent basis functions.
\item[$\bullet$] The infinite dimensional solutions of SDEs are approximated in \maxent based deterministic finite dimensional framework.
\item[$\bullet$] Error and convergence properties of \maxent based chaos expansions are studied by means of numerical simulations, and compared with the existing methods.
\item[$\bullet$] Numerical results show that if the functional dependence of a system on the random variable is unknown, \maxent based approach performs better than the polynomial expansions even if the available data is sparse.
\end{enumerate}
The rest of the paper is organized as follows. Section \ref{maxent_fcn} briefly discusses the principle of maximum entropy and the derivation of \maxent basis functions.
Solutions of SDEs using chaos expansions is reviewed in Section \ref{sec_diffEqn}. Numerical results are presented in Section \ref{sec_numerical} followed by concluding remarks in Section \ref{sec_conclude}.

\section{MAXIMUM ENTROPY BASIS FUNCTIONS} \label{maxent_fcn}
\subsection{Maximum entropy principle}
Suppose $\Delta_i$ are mutually independent events with the unknown discrete probabilities $p_i$, such that $p_i \geq 0$ and $p_i$ form the partition of unity, for $i = 1 \cdots N$.
Let us assume that we have been given expected value of a function $\E{\mb{g}(\Delta)}:=\sum_{i=1}^N p_i \mb{g}(\Delta_i)$. Our goal is to determine the probability distribution $\mb{p} := [p_1, p_2 \cdots p_N]^T$, which satisfies the given constraint.
There can be a number of distinct distributions which satisfy the given constraints.

The maximum entropy principle serves as a tool to infer the least biased distribution from the insufficient data.
To be specific, the \maxent principle states that \big(see \cite{jaynes1,jaynes2}\big), out of all possible distributions that satisfy the given constraints, the least biased or the most likely probability distribution ($\mb{p}^*$), is the one which maximizes the information-theoretic entropy $H(\mb{p})$ defined as
\begin{align}
  H(\mb{p}) := -\sum_{i=1}^N p_i \log(p_i), \eqnlabel{entropy}
\end{align}
and $0\log(0):=0$. The entropy maximization problem, formally can be written as
\begin{subequations}
\begin{equation}
  \mb{p}^* = \operatorname*{arg\,max}_{\mb{p}} H(\mb{p}),
\end{equation}
\begin{align}
  \text{such that } \sum_{i=1}^N p_i \mb{g}(\Delta_i) &= \E{\mb{g}(\Delta)} , \text{ and } \sum_{i=1}^N p_i = 1.
  \eqnlabel{subconstraintME}
\end{align}
  \eqnlabel{constraintME}
\end{subequations}
Equivalently, the constraints in \eqn{subconstraintME} can be rewritten as
\begin{align}
   \sum_{i=1}^N p_i \big( \mb{g}(\Delta_i) - \E{\mb{g}(\Delta)} \big) = 0.
  \eqnlabel{eq_subconstraintME}
\end{align}

\subsection{Minimum relative entropy principle}
Suppose prior distribution $\m:= [m_1, m_2 \cdots m_N]^T$ which estimates the probability distribution $\mb{p}$ is known, then the least biased probability distribution, $\mb{p}^*$, is determined by minimizing relative (cross) entropy or Kullback-Leibler divergence as \big(see \cite{kl_div}\big)
\begin{equation}
  \mb{p}^* = \operatorname*{arg\,min}_{\mb{p}} \Bigg( \bar{H}\big(\mb{p},\m\big):=  \sum_{i=1}^N p_i \log\Big(\frac{p_i}{m_i}\Big)   \Bigg),
  \eqnlabel{relEntropy}
\end{equation}
subject to constraints in \eqn{eq_subconstraintME}.

As mentioned earlier, \maxent basis functions that we use in the current work arise as shape functions in the context of polygonal interpolation, which is discussed next. 

\subsection{Polygonal interpolation using maxent principle}
An analogy between determining the least biased probability distribution and polygonal interpolation was first presented in \cite{sukumar}, followed by extensions of the work in  \cite{sukumarRelEntr2,ortiz}.
We briefly discuss the polygonal interpolation using \maxent basis functions below.

Let $\set{S}:=\{\param_i\}_{i=1}^N \subset \domainD \subset \Real ^{\nx}$, be a set with cardinality $N$ in $\nx-$dimensional space. Suppose that $\domainD$ is compact.
Each element $\param_i$ in $\set{S}$ is associated with a \textit{shape function} $\psi_i(\cdot) \geq 0$.
The polygonal interpolant $\hat{f}(\cdot)$ of $f(\cdot)$ at $\param$ is defined as
\begin{align}
  \hat{f}(\param) = \sum_{i=1}^N \psi_i(\param) f(\param_i), \eqnlabel{polyInterp}
\end{align}
where $f(\cdot)$ is a scalar real valued function defined over $\domainD$, and $\param \in \text{Conv}(\set{S})$, and the convex hull of the set $\set{S}$ is defined as
\begin{multline*}
\text{Conv}(\set{S}):= \{ \param | \param = \sum_{i=1}^N w_i \param_i, \sum_{i=1}^N w_i =1 , w_i \geq 0, \\ \param_i \in \set{S} \} .
\end{multline*}
The values of $\mb{\Psi}(\param):= [\psi_1(\param), \psi_2(\param) \cdots \psi_N(\param)]^T$ used in evaluation of the interpolant \eqn{polyInterp} are also referred to as barycentric coordinates of $\param$ w.r.t. elements nodes $\param_i$  of  $\set{S}$.
In this paper, we refer elements $\param_i \in \set{S}$ as basis nodes, and associated shape functions $\psi_i(\cdot)$ as basis functions.

It is desirable for interpolants to recover constant and linear functions exactly. Therefore, following constraints are imposed on barycentric coordinates $\psi_i(\param)$ to guarantee the linear precision of the interpolant.
\begin{subequations}
\begin{equation}
  \sum_{i=1}^N \psi_i(\param) = 1,
    \eqnlabel{unityPart}
\end{equation}
\begin{equation}
  \sum_{i=1}^N \psi_i(\param) \param_i = \param, \text{ or, } \Sd \mb{\Psi}(\param) = \param,
  \eqnlabel{linearConstr}
\end{equation}
\eqnlabel{linearConstrAll}
\end{subequations}
where $\Sd := [\param_1, \param_2 \cdots \param_N]$. Or equivalently,
\begin{equation}
 \tSd \mb{\Psi}(\param) = \mb{0},
  \eqnlabel{eq_linearConstr}
\end{equation}
where $\tparam_i:= \param_i-\param$ and $\tSd := [\tparam_1, \tparam_2 \cdots \tparam_N]$.

After a direct comparison between \eqn{eq_subconstraintME} and \eqn{eq_linearConstr}, $\psi_i(\param)$ can be interpreted as probability associated with the \textit{event} $\param_i$, and here $\mb{g}(\param) = \param$.
Let $\m(\param):= [m_1, m_2 \cdots m_N]^T$ be a suitably defined prior estimate for the barycentric coordinates $\mb{\Psi}(\param)$.
Therefore, $\mb{\Psi}(\param)$ which minimizes the relative entropy defined by \eqn{relEntropy} can be interpreted as the least biased barycentric coordinates of $\param$ w.r.t. basis nodes $\param_i \in \set{S}$ for a given prior estimate $\m(\param)$.

The relative entropy minimization problem is formally written as
\begin{equation}
  \mb{\Psi}^*(\param) = \operatorname*{arg\,min}_{\mb{\Psi}(\param)} \bar{H}\big(\mb{\Psi}(\param), \m(\param)\big); \text{s.t. } \tSd \mb{\Psi}(\param) = \mb{0}  \eqnlabel{PhiMEConstr_loc}
\end{equation}
Solution of this convex optimization problem using the method of Lagrange multipliers is discussed in \cite{jaynes1,sukumar,ortiz}.
This optimization problem reduces to the following system of nonlinear equations in Lagrange multipliers $\boldsymbol{\lambda}^T:=[\lambda_1, \lambda_2 \cdots \lambda_{\nx}]$ associated with equality constraints in \eqn{PhiMEConstr_loc}
\begin{align}
 \sum_{i=1}^N \tparam_i m_i(\param) e^{-\boldsymbol{\lambda}^T \tparam_i}  = \mb{0}. \eqnlabel{lamEqns_loc}
\end{align}
Numerical solvers are employed to solve this system of equations.
The optimal barycentric coordinates or basis functions in terms of $\boldsymbol{\lambda}$ are given by
\begin{align}
\psi_i^*(\param) = \frac{m_i(\param) e^{-\boldsymbol{\lambda}^T \tparam_i}}{\sum_{j=1}^N m_j(\param) e^{-\boldsymbol{\lambda}^T \tparam_j}}.
\eqnlabel{psi_loc}
\end{align}

In this paper, we use Gaussian prior, i.e. $\m(\param)$ for any $\param \in \text{Conv}(\set{S})$ is defined as
\begin{align*}
\m(\param):= [m_1, m_2 \cdots m_N]^T \text{ and } m_i(\param):= \frac{e^{-\beta ||\tparam_i||_2^2}}{\sum_j e^{-\beta ||\tparam_j||_2^2}},
\end{align*}
where $\beta \geq 0$. For $\beta > 0$, we get so called \textit{local maxent} basis functions.
The parameter $\beta$ affects the degree of locality of basis functions.
Higher value of $\beta$ implies the larger decay of basis functions away from their associated nodes, and larger degree of locality of basis functions.
For $\beta=0$ or uniform prior, i.e. $m_i = 1/N$ for $i = 1\cdots N$, we recover so called \textit{global maxent} basis functions introduced in \cite{sukumar}.
All numerical results presented in this paper are obtained for $\beta=0$.
For notational convenience, hereafter we drop the asterisk in \eqn{psi_loc} and use $\mb{\Psi}(\param)$ or $\psi_i(\param)$ to denote the optimal local \maxent basis functions evaluated at $\param$.

\section{Solution of Stochastic Differential Equations} \label{sec_diffEqn}
Let $\param\in\domainD\subset\real^d$ be a random vector with joint probability density function $\pdfp$. Approximant of a function $\vo{f}(\param):\domainD\mapsto\real^n$ using a known set of basis functions $\{\phi_i(\param)\}_{i=1}^{N}$ can be written as
\begin{align}
\vo{f}(\param) \approx \fhat(\param) :=  \sum_{i=0}^N \vo{f}_i\phi_i(\param), \eqnlabel{pcF}
\end{align}
where, $\vo{f}_i$ are deterministic coefficients. The optimal coefficients are determined by Galerkin projection, i.e. the projection of the error $\vo{e}(\param):=\vo{f}(\param) - \fhat(\param)$ against each basis function is set to zero, i.e.
\begin{align}
\E{\vo{e}(\param) \phi_i(\param)} := \int_{\domainD}\vo{e}(\param) \phi_i(\param) \pdfp d\param = \vo{0}, \eqnlabel{galerkinProj}
\end{align}
for $i=0,\cdots,N$. This results in a deterministic system of equations which can be solved for $\f_i$.
If instead of $\pdfp$, samples of $\param$ are available, then integral in \eqn{galerkinProj} can be replaced with summation and it becomes a least squares problem.

In generalized polynomial chaos (gPC) expansions, $\phi_i(\cdot)$ are multivariate orthogonal polynomials which are selected based on the known distribution of $\param$, see \cite{wiener_askey}. 
However, if $\pdfp$ is unknown, instead, samples of random vector $\param$ are available then one has to rely on data-driven methods. A data-driven polynomial chaos approach in which orthogonal basis polynomials are constructed from raw moments of available samples, termed as arbitrary polynomial chaos (aPC), is presented in \cite{aPC}.
In present work, we employ \maxent functions derived in Section \ref{maxent_fcn} as bases for chaos expansions.
Data-driven function approximation using \maxent basis functions is discussed in detail in \cite{maxentFunApprox}. In this paper, we focus on the application of \maxent basis functions for solution of stochastic differential equations.

A general procedure for solution of stochastic differential equations using chaos expansions is discussed in \cite{wiener_askey}. In this paper, for simplicity and brevity of discussion, and without loss of generality, we consider linear stochastic ordinary differential equation given by
\begin{equation}
\xdot(t,\param) = \Ap\x(t,\param),\eqnlabel{linDyn}
\end{equation}
where $\x:=\x(t,\param) \in\real^n$ and $\Ap \in \real^{n\times n}$. The solution $\x(t,\param)$ is approximated using chaos expansion as
\begin{align}
\x(t,\param) \approx \vo{\hat{x}}(t,\param) := \sum_{i=0}^N\x_i(t)\phi_i(\param) =: \X(t)\vo{\Phi}(\param),\eqnlabel{chaosExpansion} 
\end{align}
where $\X(t):=[\x_1(t) \cdots \x_N(t)]$ are time varying deterministic coefficients and $\vo{\Phi}(\param):= [\phi_1(\param) \cdots \phi_N(\param)]^T$ are chosen basis functions for chaos expansion.
After substituting the approximation \eqn{chaosExpansion}, the projection of equation error of \eqn{linDyn} against each basis is set to zero to obtain a system of deterministic ODEs in terms of chaos coefficients $\x_i(t)$, which is given by
\begin{multline}
\xdot_c = \left(\E{\phip\phipt}\otimes\I{n}\right)^{-1} \times \\ \E{(\phip\phipt)\otimes\Ap)}\x_c,
\eqnlabel{gpcLinSurr}
\end{multline}
where, $\x_c(t):= \vect{\X(t)}$.
A step-by-step derivation of the surrogate system of deterministic ODEs similar to \eqn{gpcLinSurr} for a generalized non-linear ODE can be found in \cite{accuratePC}.

The deterministic ODEs in \eqn{gpcLinSurr} can be solved using standard deterministic approaches. Therefore, $\x_c(t)$ and hence $\X(t)$ can be determined at any time $t$. Then the moments of $\x(t,\param)$ can be approximated by the moments of $\vo{\hat{x}}(t,\param)$ as
\begin{subequations}
\begin{equation}
\E{\x(t)} \approx \E{\vo{\hat{x}}(t)} = \X(t)\E{\vo{\Phi}(\param)},\\
\end{equation}
\begin{multline}
\mathbb{E}\big[ {\big(\x(t) - \E{\x(t)}\big) \big(\cdot\big)^T} \big]  \approx \mathbb{E}\big[ {\big(\hat{\x}(t) - \E{\hat{\x}(t)}\big) \big(\cdot\big)^T} \big] \\
= \X(t) \Big(\E{\phip\phipt} -\E{\phip}\E{\phip}^T \Big) \X^T(t).
\end{multline}
\eqnlabel{approxMoments}
\end{subequations}

Let $\set{D}:=\{\param_i\}_{i=1}^{\nd} \subset \domainD$ be the set of available samples for $\param$. Since this is the only available information about $\param$, the expectation integrals involved in \eqn{gpcLinSurr} and \eqn{approxMoments} are approximated by sample averages calculated using data points in $\set{D}$.

We define $\set{B}:= \{\param_i\}_{i=1}^{\nb}$ such that $\set{D}\subseteq \cvxh{\set{B}}$. This is trivially satisfied if we choose $\set{B} = \set{D}$.
However, $\set{B}$ need not be identical to $\set{D}$. The elements in $\set{B}$ serve as nodes for \maxent basis functions. Therefore, the number of \maxent basis functions is $\nb$. The definition of $\set{B}$ ensures that barycentric coordinates (see \eqn{psi_loc}) of each $\param_i \in \set{D}$ are well defined w.r.t. basis nodes in $\set{B}$.

In the following section we solve different stochastic ODEs using \maxent basis functions as bases for chaos expansions.

\section{Numerical results} \label{sec_numerical}

For the purpose of numerical simulations, we consider a scalar ODE given by
\begin{align}
\dot{x}(t) = a(\Delta)x(t), \quad \Delta \in \real, \eqnlabel{exODE}
\end{align}
with deterministic initial condition $x(0) = 1$. The decay coefficient $a(\Delta)$ is a stochastic parameter with dependence on random variable $\Delta$.
The system of ODEs for chaos coefficients when using \maxent basis functions follows from \eqn{chaosExpansion} and \eqn{gpcLinSurr} as
\begin{align}
&x(t,\Delta) \approx \hat{x}(t,\Delta) = \sum_{i=0}^{\nb} x_i(t)\psi_i(\Delta), \nonumber\\
&\xdot_c = \left(\E{\mb{\Psi}(\Delta)\mb{\Psi}^T(\Delta)}\right)^{-1}\E{a(\Delta)\mb{\Psi}(\Delta)\mb{\Psi}^T(\Delta)}\x_c, \eqnlabel{gpcScalarLinSurr}
\end{align}
where $\x_c(t) = [x_1(t), x_2(t) \cdots x_{\nb}(t)]^T$. Note, $\mb{\Psi}(\Delta)$ are barycentric coordinate of $\Delta$ w.r.t basis nodes in $\set{B}$. Approximate mean and variance at time $t$ follow from \eqn{approxMoments} as
\begin{subequations}
\begin{equation}
\hat{\mu}_x(t) = \x_c(t)^T\E{\mb{\Psi}(\Delta)},
\end{equation}
\begin{equation}
 \hat{\sigma}_x^2(t) = \x_c^T(t)(\E{\mb{\Psi}(\Delta)\mb{\Psi}^T(\Delta)} -  \E{\mb{\Psi}(\Delta)}\E{\mb{\Psi}(\Delta)}^T ) \x_c(t)
\end{equation}
\eqnlabel{scalarApproxMom}
\end{subequations}
Since \maxent basis functions evaluated for any $\Delta$ add up to one (see \eqn{unityPart}), the initial condition for chaos coefficients can be written as $\x_c(0) = \mb{1}_{\nb}$, where, $\mb{1}_{\nb}$ is the column vector of length $\nb$ with each element unity.
With this initial condition for $\x_c$, it is straightforward to verify that mean and variance approximated using \eqn{scalarApproxMom} come out to be $\hat{\mu}_x(0) = 1$ and $\hat{\sigma}_x^2(0) = 0$, which exactly recovers the given initial condition $x(0)=1$.

In text to follow, we consider two examples. The first example is a classical ODE studied in chaos expansion literature. In this example, we assume that the functional form of $a(\Delta)$ is known and investigate the error convergence properties of \maxent based expansions.
In second example, we assume that the functional form of $a(\Delta)$ is unknown, instead, an approximant is constructed from a given sparse data set. Then a surrogate model is developed using the approximant. The results obtained for \maxent based models are compared with data-driven polynomial chaos (aPC) based surrogate models.

\subsection*{Example 1}
Let us assume that the functional form of $a(\Delta) = -(1+\Delta)/2$ is known, where $\Delta$ is uniformly distributed over $[-1,1]$. 
This ODE has been used as a test problem in \cite{wiener_askey, aPC}.
The analytical or true values of mean and variance are given by $\mu_x(t) = \frac{1-e^{-t}}{t}$ and $\sigma_x^2(t) = \frac{1-e^{-2t}}{2t}-\left(\frac{1-e^{-t}}{t}\right)^2$ for $t>0$,
which are shown by black dashed lines in \fig{unif_mean} and \fig{unif_var}.

The two major sources of error in estimated moments given by \eqn{scalarApproxMom} are namely, limited number of data samples ($\nd$), and truncation of chaos expansion at finite number of terms or finite number of basis functions ($\nb$). Convergence properties of both types of error are discussed in the following text.
We first consider the effect the of number of basis functions, and for this investigation we assume the data set $\set{D}$ to be fixed. For the second case, we keep the basis functions fixed, and vary the size of data set $\set{D}$.

\fig{unif_mean} and \fig{unif_var} show the effect of increasing number of \maxent basis functions for a fixed data set $\set{D}$ which consists of $\nd=500$ uniformly spaced points in $[-1, 1]$. The set of basis nodes, $\set{B}$, is selected as a collection of $\nb$ uniformly spaced points in $[-1,1]$.
\Eqn{gpcScalarLinSurr} is integrated numerically from $t=0$ to $t=30$ for different values of $\nb$. Moments at each time instant are approximated using \eqn{scalarApproxMom} and are shown in plots on the left sides of \fig{unif_mean} and \fig{unif_var}.
It is well known that moments approximated using surrogate model such as \eqn{gpcScalarLinSurr} deviate from true values as ODEs are integrated forward in time. The errors in estimated moments calculated w.r.t. true moments are shown in plots on the right sides of \fig{unif_mean} and \fig{unif_var}. The error in estimated mean is defined as $\varepsilon_{\mu}(t):= |1- \hat{\mu}_x(t) /\mu_x(t)|$.
The error in estimated variance is defined in similar way. It is clear from \fig{unif_mean} and \fig{unif_var} that surrogate models obtained using more number of \maxent basis functions maintain better accuracy for longer periods of time while performing the temporal integration.

\begin{figure}
\begin{center}
\includegraphics[trim={1cm 0cm 1.3cm 0cm},clip,width=0.48\textwidth]{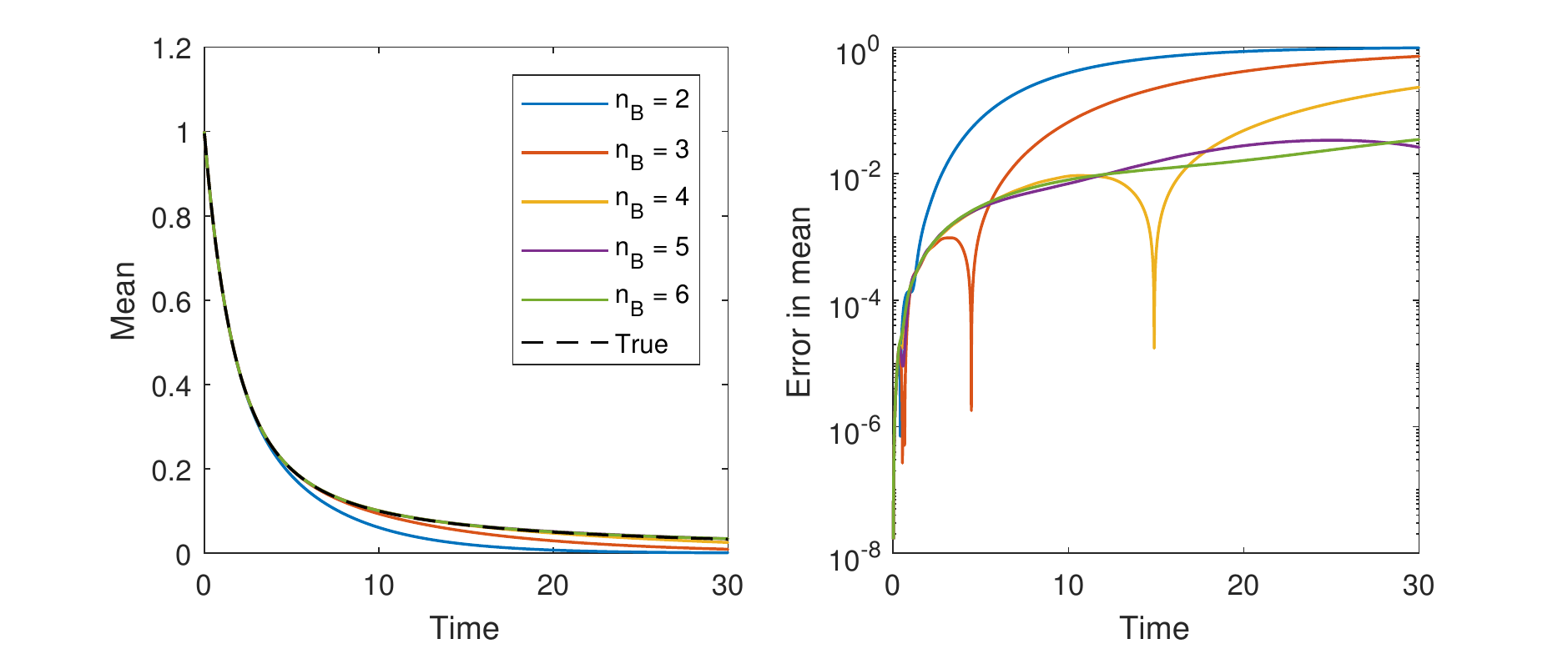}    
\caption{Estimated mean and error for different number of \maxent basis functions ($\nb$) and fixed $\nd=500$.}
\figlabel{unif_mean}
\end{center}
\end{figure}

\begin{figure}
\begin{center}
\includegraphics[trim={1cm 0cm 1.3cm 0cm},clip,width=0.48\textwidth]{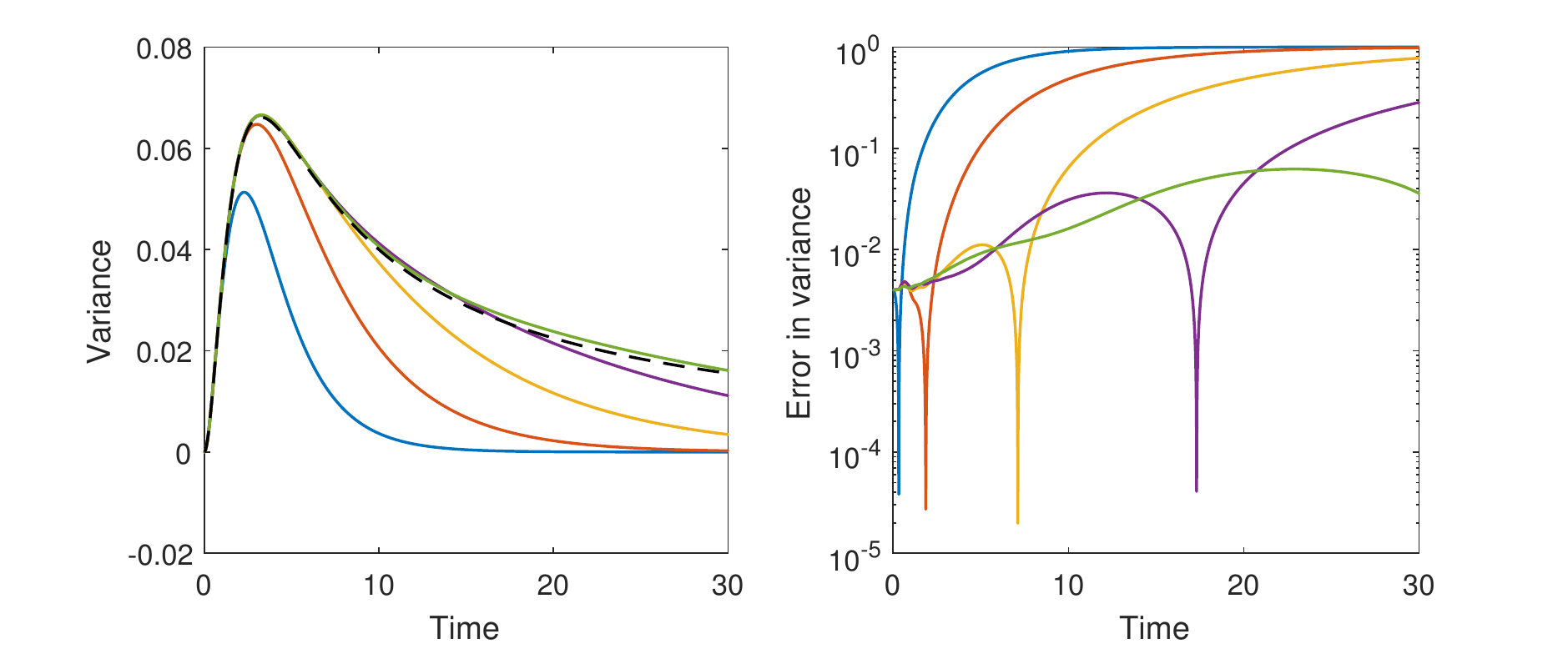}    
\caption{Estimated variance and error for different number of \maxent basis functions ($\nb$) and fixed $\nd=500$.}
\figlabel{unif_var}
\end{center}
\end{figure}

\fig{err_conv} shows the convergence of error in moments calculated at $t=10$ w.r.t. number of basis functions ($\nb$). It is clear that both \maxent and aPC based models demonstrate similar error convergence trend.
We observe strong convergence initially as we start increasing number of basis functions. However, after certain value of $\nb$, the errors saturate and do not decreases further as $\nb$ is increased.
For example, consider error in mean shown by blue line with circles in \fig{err_conv}.
We observe strong convergence as the number of basis functions is increased from $\nb=2$ to $\nb=4$. However, for $\nb=5$ onwards, the errors become stagnant.
This observation is attributed to limited size of sample set. This suggests that given a finite set of samples, $\set{D}$, there is a lower limit on how much error can be reduced.

\begin{figure}
\begin{center}
\includegraphics[trim={0.2cm 0cm 0.9cm 0cm},clip,width=0.27\textwidth]{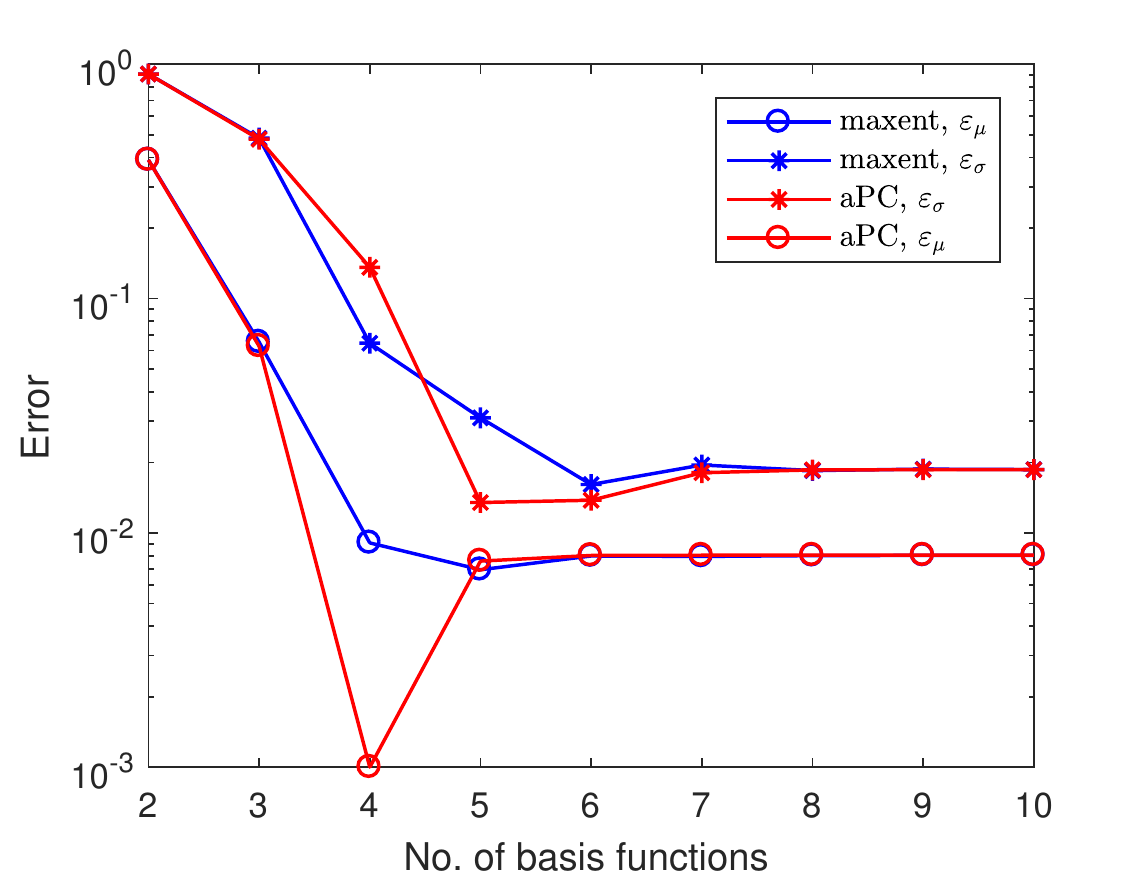}    
\caption{Error in estimated moments at $t=10$ for different number of basis functions ($\nb$) and fixed $\nd=500$.}
\figlabel{err_conv}
\end{center}
\end{figure}

\fig{compare_apc} compares temporal evolution of error in estimated moments for data-driven surrogate models obtained using aPC and the proposed \maxent based framework, for equal number of basis functions.
We observe that \maxent and aPC based models demonstrate similar error evolution, with \maxent based model being marginally better than the aPC based model for estimated variance for $t>10$.

\begin{figure}
\begin{center}
\includegraphics[trim={0.7cm 0cm 1.3cm 0cm},clip,width=0.48\textwidth]{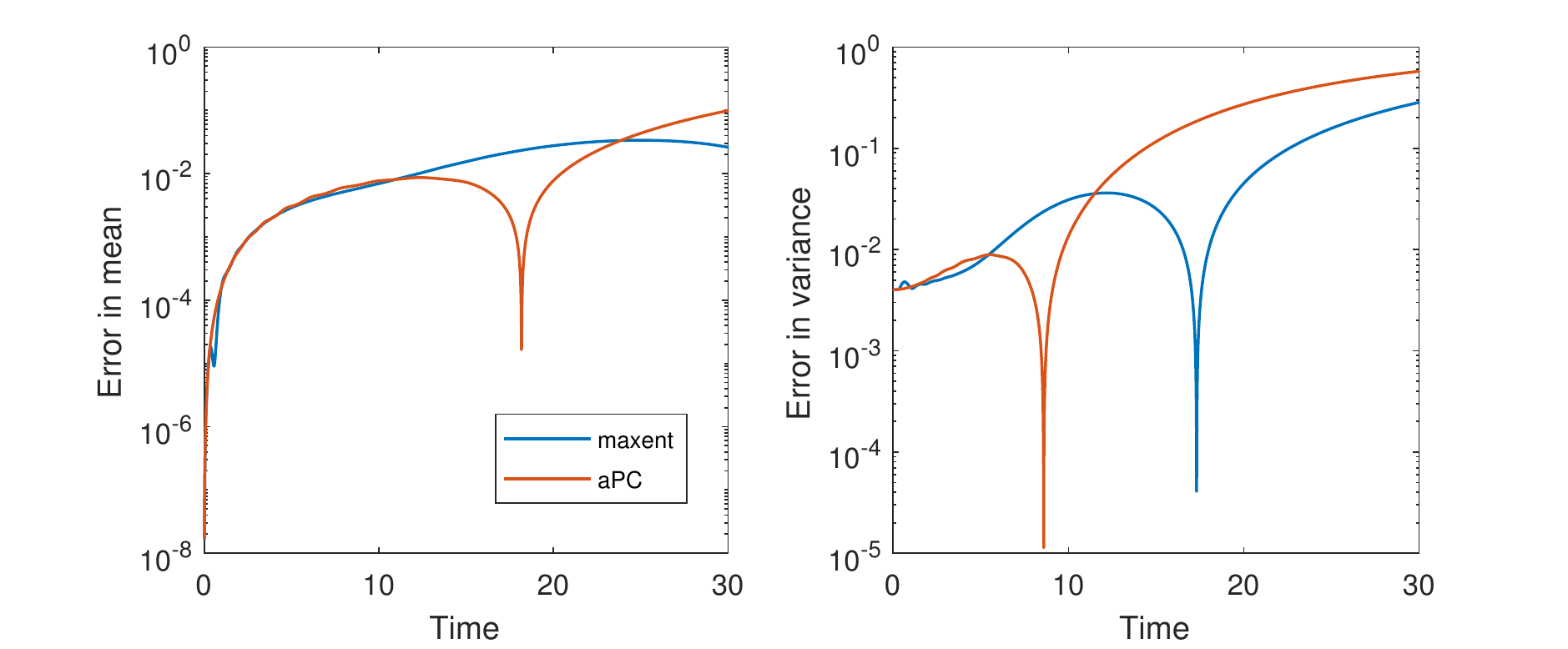}    
\caption{Comparison of error in estimated moments for aPC and \maxent based models, $\nb = 5$ and $\nd = 500$.}
\figlabel{compare_apc}
\end{center}
\end{figure}

Needless to say, the accuracy of approximation depends on the number of available samples, i.e. $\nd$.
The left plot in \fig{err_conv_samp} shows the mean of errors in estimated mean and variance of $x(t)$ at $t=10$ calculated for different number of data points ($\nd$) sampled randomly from $[-1,1]$, repeated 500 times, for a fixed number of \maxent basis functions $\nb=5$.
Similar variation for the variance of errors in estimated moments is shown in the right plot.
With no surprise, we observe that both mean and variance of errors decrease as the number of available samples increases.
The aPC based model also shows a very similar trend (not shown in the figure).

\begin{figure}
\begin{center}
\includegraphics[trim={0.7cm 0cm 1.3cm 0cm},clip,width=0.48\textwidth]{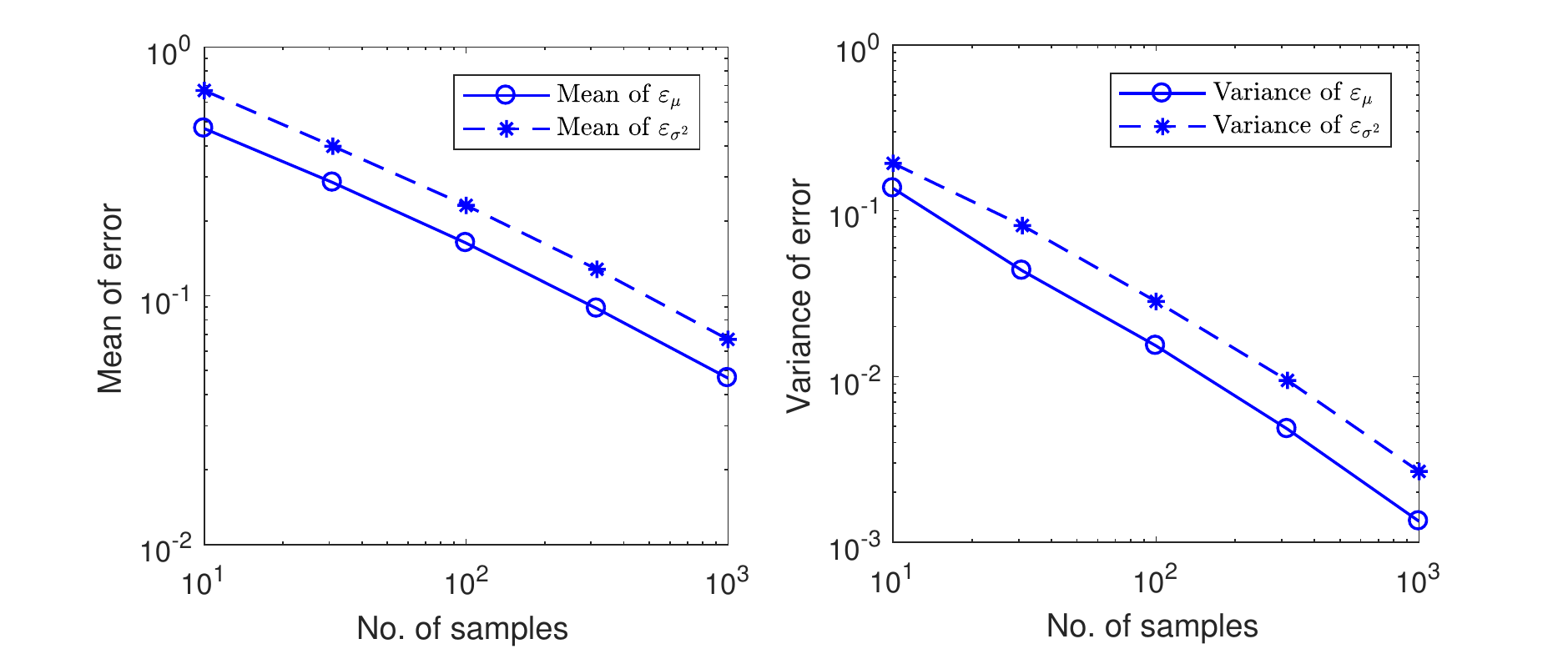}    
\caption{Variation of mean (left) and variance (right) of error in estimated moments with different number of samples ($\nd$), for fixed $\nb = 5$ \maxent basis functions, at $t=10$.}
\figlabel{err_conv_samp}
\end{center}
\end{figure}

Results discussed herein are obtained if $\Delta$ is sampled from a uniform distribution. Although not shown here, similar results are observed when $\Delta$ is sampled from a Gaussian distribution.
In this example, where functional form of $a(\Delta)$ is assumed to be known, we observe that aPC and \maxent based models demonstrate similar performances in terms of accuracy of the estimated moments and error convergence rates.

\subsection*{Example 2}
In this example we assume that the functional form of $a(\Delta)$ is not known. Instead, we have been given samples of $\Delta$, and values of $a(\Delta)$ for a sparse sample set of $\Delta$. We assume that $\nd$ samples of $\Delta$ are given.
Let $\set{D}':= \{\Delta_j\}_{j=1}^{n_{\set{D}'}}$ be the sparse set of samples of $\Delta$ for which $a_j:=a(\Delta_j)$ is known.

Here, solving stochastic ODE involves two steps. First,  we approximate the function $a(\Delta)$ using given data $\{\Delta_j,a_j\}_{j=1}^{n_{\set{D}'}}$. Second, we construct chaos expansion as we did in the previous example and integrate the surrogate model temporally.
Basis nodes for \maxent basis functions are chosen to be the elements in $\set{D}'$, i.e. $\set{B}$ and $\set{D}'$ are identical and $\nb = n_{\set{D}'}$. The same basis nodes are used for function approximation in step one and chaos expansion in step two.
The weights or coefficients associated with basis functions for function approximation are determined by the least squares solution as discussed in \cite{maxentFunApprox}. 

For the purpose of simulation, we assume that $\Delta$ is a uniformly distributed variable between $(0,1)$, and that $\set{D}$ and $\set{D}'$ respectively consist of uniformly spaced $\nd=500$ and $n_{\set{D}'} = 10$ points in $(0,1)$. The unknown underlying functional form of $a(\Delta)$ is assumed to be $a(\Delta) = \Delta^{\alpha-1}(1-\Delta)^{\gamma-1}$ with $\alpha = 2, \gamma = 0.5$.
Moments of the true or reference solution of \eqn{exODE} are obtained using Monte-Carlo runs with $5\times 10^4$ samples of $\Delta$.

\begin{figure}
\begin{center}
\includegraphics[trim={0.7cm 0cm 1.3cm 0cm},clip,width=0.48\textwidth]{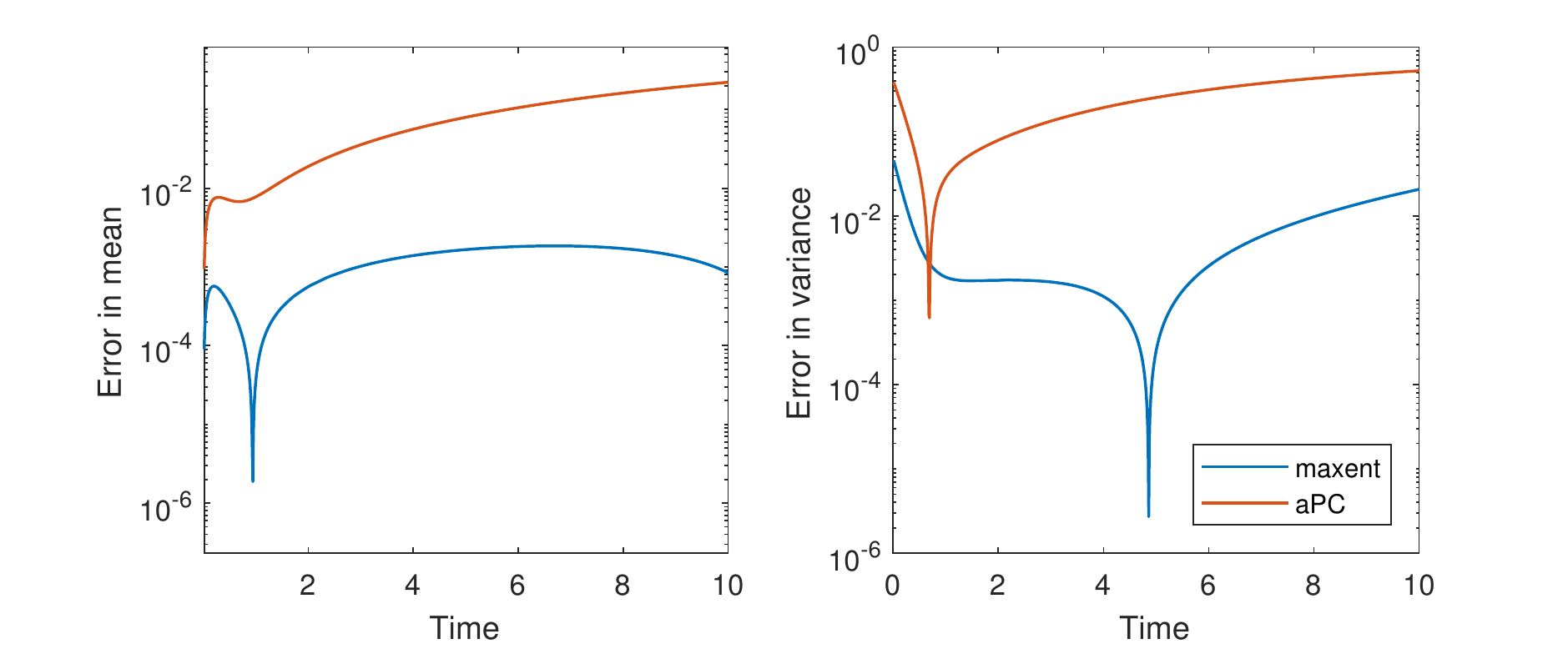}    
\caption{Comparison of error in estimated moments, $\nb = n_{\set{D}'} = 10$, $\nd = 500$.}
\figlabel{param_fcnApp}
\end{center}
\end{figure}

Comparison of error in estimated moments using \maxent and aPC based models is shown in \fig{param_fcnApp}, and clearly, the former demonstrates better accuracy. The relatively better accuracy of the \maxent model is mainly attributed to the ability of \maxent basis functions to approximate the underlying functional form of $a(\Delta)$ with sparse data
better than polynomials bases constructed for the aPC model.

The normalized error for function approximation is shown in \fig{err_fcnApp}. The \maxent based approximant has lower error than the aPC based function approximant, and consequently, \maxent based model demonstrates better accuracy in the estimated moments as shown in \fig{param_fcnApp}.
However, as discussed in \cite{maxentFunApprox}, it should be noted that if the underlying functional form of $a(\Delta)$ lies within the span of basis polynomials, then aPC based model will have similar, or perhaps even better accuracy than the \maxent model.
But, for functions such as we considered here, which are difficult to approximate using polynomials, especially if the given data set is sparse, \maxent based models perform better.
Ability to accurately approximate functions is particularly vital for solving differential equations with the initial condition uncertainty, as modeling of uncertainty in the initial condition reduces to a function approximation problem.

\begin{figure}
\begin{center}
\includegraphics[trim={0.2cm 0cm 1.1cm 0cm},clip,width=0.27\textwidth]{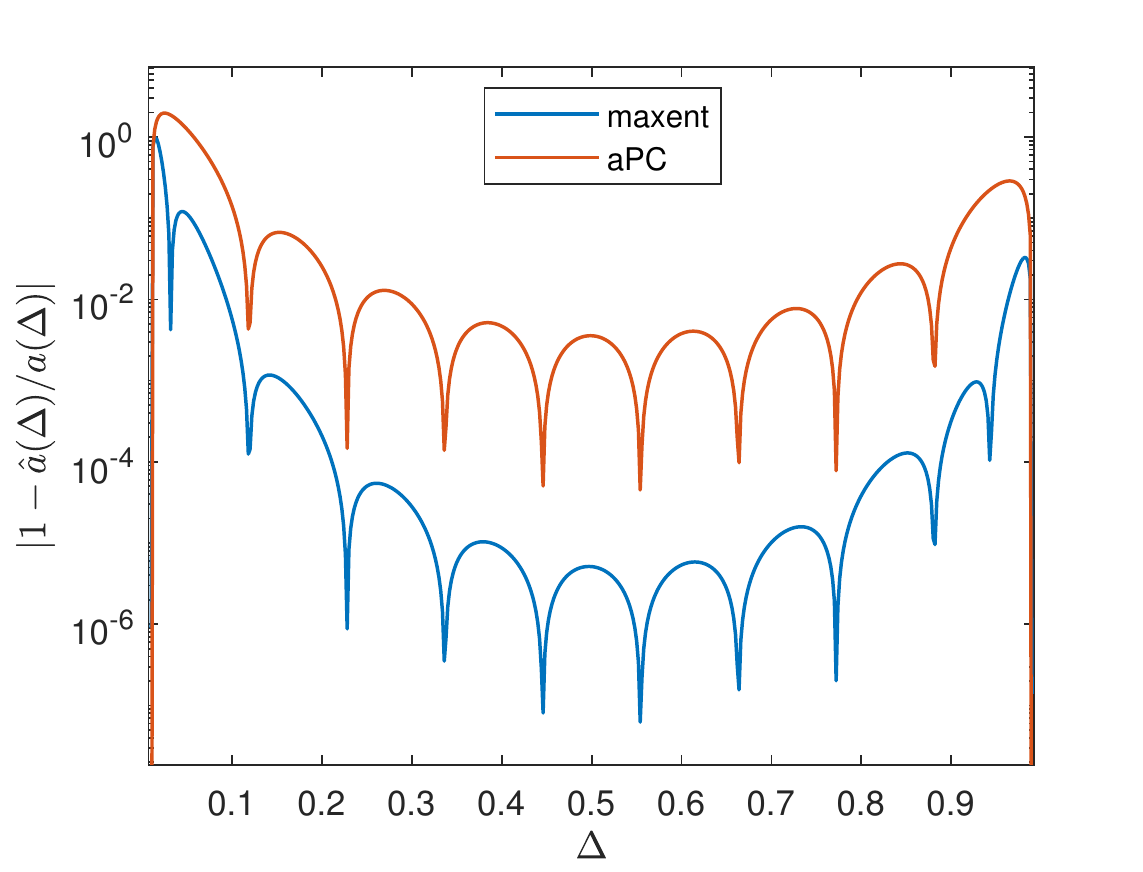}    
\caption{Comparison of function approximation accuracy, $a(\Delta) = \Delta/\sqrt{1-\Delta}$, $\nb = n_{\set{D}'} = 10$.}
\figlabel{err_fcnApp}
\end{center}
\end{figure}

\section{Conclusion} \label{sec_conclude}
In this paper we presented an approach to develop surrogate models based on limited data by construction of chaos expansions using \maxent basis functions.
We investigated the error characteristics and convergence properties of such \maxent based chaos expansions and compared the accuracy with data-driven polynomial chaos expansions (aPC).
We observed that \maxent based surrogate models demonstrate similar accuracy as the aPC models if the functional dependence on the random variables is known.
In the case where functional dependence is unknown, \maxent based approach demonstrated better accuracy than the polynomial expansions especially if the available data is sparse.

We note that the selection of basis nodes for \maxent basis functions may not be trivial for high dimensional complex systems, and it may affect the accuracy of chaos expansion. However, this investigation is out of scope of this paper. 
In this paper, we considered simple differential equations, with data set sampled from uniform and Gaussian distributions.
However, results observed in this paper are motivating enough to pursue further investigations with more complex underlying distributions and differential equations of real physical systems, and will be topics of our future works.
\vspace{-0.2cm}
\bibliographystyle{plainurl}
\bibliography{root}             


\end{document}